\documentclass[11pt]{article}
\usepackage{amssymb}
\usepackage{amsmath}\usepackage{graphicx}
\usepackage{times}\hfuzz=10pt 
\newtheorem{theorem}{Theorem}%[section]
\newtheorem{lemma}{Lemma}%[section]
\newtheorem{definition}{Definition}%[section]
\newtheorem{remark}{Remark}%[section]
\def\e{\varepsilon}

\def\defi{\stackrel{{\scriptscriptstyle \Delta}}{=}}

\def\o{\omega}
\def\O{\Omega}

\def\Y{{\cal Y}}

\def\w{\widehat}
\def\Ind{{\mathbb{I}}}

\def\esssup{\mathop{\rm ess\, sup}}

\def\Im{{\rm Im\,}}

\def\R{{\bf R}}

\def\S{{\bf S}}
\def\Z{{\cal Z}}
\def\ZZ{{\bf Z}}

\def\g{\gamma}
\def\C{{\bf C}}
\def\W{{\cal W}^*}
\def\W{{\cal W}}
\def\ww{\widetilde}
\def\X{{\cal X}}
\def\t{\theta}
\def\oo{\bar}

\def\U{{\cal U}}
\def\V{{\cal V}}

\def\T{{\mathbb{T}}}

\newcommand{\be}{\begin{equation}}
\newcommand{\ee}{\end{equation}}
\newcommand{\bd}{\begin{displaymath}}
\newcommand{\ed}{\end{displaymath}}
\newcommand{\ba}{\begin{array}{ll}}
\newcommand{\ea}{\end{array}}
\newcommand{\baa}{\begin{eqnarray}}
\newcommand{\eaa}{\end{eqnarray}}
\newcommand{\baaa}{\begin{eqnarray*}}
\newcommand{\eaaa}{\end{eqnarray*}}

\def\oo{\bar}

\def\K{{\cal K}}
\def\ew{\left(e^{i\o}\right)}
\date{ Submitted: August 25, 2014. Revised: September 23, 2014 }% xxx
\title{
Causal band-limitness and predictability criterions for one-sided sequences}
\author{
Nikolai Dokuchaev}
 
\begin{document}
 \vspace{-0.5cm}      \maketitle
\def\brea{}
\def\breakk{}
\def\break{}
\begin{abstract} The paper studies frequency characteristics and predictability of  real sequences, i.e.,  discrete time processes in deterministic setting. We consider  band-limitness and predictability of one-sided sequences.
We establish predictability of some classes of sequences, and we suggest  predictors represented by
causal convolution sums over past times.
\index{in
IEEE:}\let\thefootnote\relax\footnote{The author is with  Department of
Mathematics and Statistics, Curtin University, GPO Box U1987, Perth,
Western Australia, 6845 (email N.Dokuchaev@curtin.edu.au).  This work  was
supported by ARC grant of Australia DP120100928 to the author.}
\par MSC 2010 classification : 42A38, %       Fourier and Fourier-Stieltjes transforms and other transforms of Fourier type
93E10, %estimation and detection
62M15,      % Spectral analysis is stochastic processes
42B30  %Hardy spaces
\end{abstract}\section{Introduction}
The paper studies frequency characteristics and predictability of  real sequences, i.e.,  discrete time processes in deterministic setting.
It is well known that certain restrictions on the frequency distribution can ensure some opportunities for prediction and interpolation of the processes.      For continuous time processes, the classical result  is the Nyquist-Shannon-Kotelnikov interpolation theorem
for the sampling of continuous
time band-limited processes.  These processes are analytic; they allow a unique extrapolation from any interval and are uniquely defined by their past.   A similar result was obtained for the processes  with the exponential decay of energy on the higher frequencies that are not necessary band-limited \cite{D10}.
Predictability based on sampling and the
Nyquist-Shannon-Kotelnikov theorem was discussed, e.g., in
\cite{F,K,L,Ly,M,Mi,N,NT}. These works
considered  continuous-time band-limited
stochastic processes; the corresponding
predictors were constructed for the setting
where the shape of the spectral representation is supposed to be
known. For discrete time processes, the predictability can also be
achieved given some properties of spectral representations.

For discrete time processes or sequences,  it is not obvious how to
define an analog of the continuous time analyticity. So far, there
is a criterion of predictability in the frequency domain setting
given by the classical Szeg\"o-Kolmogorov Theorem for stochastic
Gaussian stationary discrete time processes. This theorem says that
the optimal prediction error is zero  if the spectral density $\phi$
is such that \baa \int_{-\pi}^\pi \log\phi\ew d\o=-\infty; \eaa see \cite{K41,Sz,V}\index{
Szeg\"o (1920, 1921), Verblunsky  (1936), Kolmogorov (1941),\index{
Theorem 2}} and recent literature reviews in \cite{B,S}.\index{ Simon (2011) and Bingham
(2012).} This means that a stationary Gaussian  process  is
predictable if its spectral density is vanishing on a part of the
unit circle $\T=\{z\in\C:\ |z|=1\}$, i.e., if the process is
band-limited in this sense. This result was expanded on more
general stable stochastic processes  allowing spectral representations with
spectral density via processes with independent increments (see, e.g., \cite{CS}).\index{ Cambanis and Soltani (1984))}
In \cite{D10,D12a,D12b,Df,D13}, this problem was
readdressed in the deterministic setting, and some predictors
were suggested for band-limited sequences, i.e., for sequences for
with Z-transform vanishing on a part of $\T$. However, it appears that the framework of two-sided sequences
 required for
detecting of the band-limitness via Z-transforms are  not always convenient to use. For example, consider a situation where the data is collected
dynamically during a prolonged time interval. Application of the two-sided Z-transform requires to to select some
past time at the middle of the time interval of the observations as the zero point for a model of  the two-sided sequence; this could be inconvenient. Therefore, it could be more  more convenient to analyze one-sided sequences rather than  two-sided sequences required for
detecting of the band-limitness via Z-transforms.
This leads to the analysis of the  one-sided sequences directed backward to the past. However,
the straightforward application of the one-sided
Z-transform to the historical data considered as an one-sided sequence directed to the past does not generate Z-transform vanihing on a part of the unit circle even for a band-limited underlying  sequence.  To overcome this, we use  sine  and cosine versions of  Z-transforms;  it appears
 that they allow to detect band-limitness in one-sided sequences.
This gives a possibility  to  establish predictability of certain classes of one-sided sequences.
Following \cite{D10,D12a,D12b}, we suggest  predictors represented by
causal convolution sums over past times representing historical
observations.
\section{Definitions and the main results}\label{secMain}
We denote by $\ZZ$ the set of all integers.
\par
For $\tau\in\ZZ\cup\{+\infty\}$ and $\t<\tau$,
we denote by $\ell_r(\t,\tau)$ a Banach space  of
sequences $x=\{x(t)\}_{\t-1<t<\tau+1}\subset\R$, with the norm
$\|x\|_{\ell_r(\t,\tau)}=\left(\sum_{t=\t}^{\tau}|x(t)|^r\right)^{1/r}<+\infty$
for $r\in[1,\infty)$ or  $\|x\|_{\ell_\infty(\t,\tau)}=\sup_{t:\ \t-1<t<\tau+1}|x(t)|<+\infty$
for $r=+\infty$; the cases where
$\t=-\infty$ or $\tau=+\infty$ are not excluded. As usual, we assume that all sequences with the finite
 norm  of this kind are included in the corresponding space.
\par
\def\elm{\ell_2^-}
\def\elp{\ell_2^+}
\def\elr{\ell_r^-}
\def\elrr{\ell_r^{++}}
\def\CC{\Xi}
\par
For brevity, we will use the notations $\ell_r=\ell_r(-\infty,\infty)$, and $\ell_r^-=\ell_r(-\infty,0)$.\index{, $\ell_r^+=\ell_r(0,\infty)$, $\ell_r^{++}=\ell_r(1,\infty)$.}
\par
\subsection{The classical Z-transform and band-limitness}
For  $x\in \ell_1$ or $x\in \ell_2$, we denote by $X=\Z x$ its
Z-transform defined as \baaa X(z)=\sum_{t=-\infty}^{\infty}x(t)z^{-t},\quad
z\in\C. \eaaa Respectively, the inverse $x=\Z^{-1}X$ is defined as
\baaa x(t)=\frac{1}{2\pi}\int_{-\pi}^\pi X\left(e^{i\o}\right)
e^{i\o t}d\o, \quad t=0,\pm 1,\pm 2,....\eaaa
\par
Let   $\T\defi\{z\in\C:\ |z|=1\}$, and let  $\T^+\defi\{z\in\C:\ |z|=1,\quad \Im z\ge 0\}$. If $x\in \ell_2$, then $X|_\T$ is defined as an element of
$L_2(\T)$. Let
the mapping $\Z^+:\ell_2\to L_2(\T^+)$ be  defined as $\Z^+x=(\Z x)|_{\T^+}$.
It can be noted that $\overline{X\ew} =X\left(e^{-i\o}\right)$ for $X=\Z x$, $x\in\ell_2$, and the mapping $\Z^+:\ell_2\to L_2(\T^+)$ is a continuous bijection between $\ell_2$ and $L_2(\T^+)$.
\begin{definition}
We will call a two-sided sequence  $x\in\ell_2$   band-limited
if there exists $\O\in[0,\pi)$ such that $X\ew =0$ for $\o\in[-\pi,\O)\cup(\O,\pi]$, where $X=\Z x$.
\end{definition}
It is known that two-sided band-limited infinite sequences are predictable in a certain sense  (see Theorem 1 in \cite{D12a} and Theorem 1 in \cite{D12b}).

As was mentioned in Section 1, for many practical applications, it is  inconvenient to estimate frequency characteristics of two-sided infinite sequences.
On the other hand,  the unilateral Z-transform  calculated for  a one-sided part of a
band-limited process does not show its band-limitness.
More precisely, if we consider a band-limited sequence $x\in\ell_2$, then
the process $x_\t(t)=x(t)\Ind_{(-\infty,t]}(t)$ cannot be band-limited for a given $\t\in\ZZ$, and $\ln|(\Z x_\t)\ew|\in L_1(-\pi,\pi)$, there are constraints on the degeneracy of $\Z x_\t$ on $\T$.
\subsection{Sine and cosine Z-transforms for one-sided sequences}
We suggest below  sine and cosine modifications of Z-transform oriented on applications to the one-sided sequences.
\par
For  $x\in \ell_1^-$ or $x\in \ell_2^-$, we introduce transforms $\xi_1=\CC_1 x\in L_2([0,\pi],\R)$ and $\xi_2=\CC_2 x\in L_2([0,\pi],\R)\times \R$ \index{$\xi_k=\CC_k x\in L_2([0,\pi],\R)$, $k=1,2$,} defined as
 \baaa &&\xi_1(\o)=2\sum_{t=-\infty}^{-1}\cos(\o t)x(t)+x(0),\quad
\\ &&\xi_2(\o)=\left(2\sum_{t=-\infty}^{-1}\sin(-\o t)x(t),\,x(0)\right),\eaaa where $\o\in[0,\pi]$.
\begin{lemma}\label{lemma1} The mappings $\CC_1:\elm\to L_2([0,\pi],\R)$ and $\CC_2:\elm\to L_2([0,\pi],\R)\times \R$, \index{$\CC_k:\elm\to L_2([0,\pi],\R)$, $k=1,2$,} are continuous bijections, and
the corresponding inverse mappings  are also a continuous bijections such that $x_k=\Xi_k^{-1}\xi_k$ are defined as
\baa &&x_1(t)=\frac{1}{\pi}\int_{0}^\pi \xi_1\left(\o\right)\cos(\o t)d\o,\breakk\quad t=0,-1,-2,...,\nonumber\\
&&x_2(t)=\frac{1}{\pi}\int_{0}^\pi \xi_2'\left(\o\right)\sin(-\o t)d\o,\breakk\quad t=-1,-2,...\nonumber\\
&&x_2(0)=\xi_2'',\quad\hbox{where}\quad \xi_2=(\xi_2'(\o),\xi_2'').\hphantom{xx}\label{invinv}\eaa
\end{lemma}
\subsubsection{Application to band-limitness and predictability}
To justify the introduction of the new transforms $\Xi_k$,
we demonstrate below that, similarly to the  band-limitness defined via the  Z-transform for two-sided infinite sequences,
it is possible to define a detectable analog of the band-limitness for the one-sided infinite sequences. Moreover,
we will show that the presence of this new band-limitness also leads to a predictability.
\begin{definition}\label{defBL} Let $r\in\{1,2\}$.
We will call a one-sided sequence  $x\in\ell_r^-$   causally band-limited (or left band-limited)
if there exists $\O\in[0,\pi)$ and $x'\in\ell_r(1,+\infty)$  such that $X\ew =0$ for $|\o|>\O$, where $X=\Z \oo x$, and where $\oo x\in \ell_r$ is such that $\oo x(t)=x(t)$ for $t\le 0$ and $\oo x(t)=x'(t)$ for $t> 0$.
\end{definition}
In principle, it is possible to verify if $x$ the conditions of this definition holds with a given $\O\in[0,\pi)$. In particular, it can be done
using Theorem 1 from \cite{Df,D13} via solution of a infinity dimensional quadratic optimization problem which is a computationally challenging procedure. Moreover, this can be done for each potentially acceptable   $\O$ separately.

We suggest below a more convenient sufficient conditions of band-limitness.
\begin{theorem}\label{Th0} A sequence $x\in\elm$ is causally band-limited
if there exists $\O\in[0,\pi)$ such that at least one of the following two condition holds:
\begin{itemize}
\item[(i)] There exists  $a\in\R$ such that  $\xi_1(\o) =a$ for $\o\in(\O,\pi]$, where $\xi_1=\CC_1 x$;
\item[(ii)]  $\xi_2'(\o) =0$ for $\o\in(\O,\pi]$, where $\xi_2=(\xi_2'(\o),\xi_2'')=\CC_2 x$.
\end{itemize}
\end{theorem}
    It can be noted that  a sequence can be causally band-limited even if the transform $\xi_k(\o)$ are separated from zero on  $[0,\pi]$; this feature makes causal band-limitness different
from the  band-limitness defined for the two-sided sequences via degeneracy of the Z-transform.
\begin{remark}\label{rem1}
The  conditions on $\xi_1$ and $\xi_2'$ required in Theorem \ref{Th0} cannot be satisfied simultaneously.
\end{remark}
\section{Predicability of causally band-limited processes and more general processes}
Let $D\defi\{z\in\C: |z|\le 1\}$, $D^c=\C\backslash D$.\index{ $\ZZ$ is the set of all
integers} For $r\in[1,+\infty]$, let  $H^r$ be the Hardy space of functions that are holomorphic on
$D^c$ including the point at infinity   (see, e.g.,  \cite{Du}). Note that Z-transform defines a bijection
between the sequences from $\ell_2^+$ and the restrictions (i.e.,
traces) $X|_{\T}$ of the functions from $H^2$ on $\T$ such that  $\overline{X\ew} =X\left(e^{-i\o}\right)$.

\begin{definition}
Let $\w\K$ be the class of functions $\w k:\ell_{\infty}^+\to\R$
such that $\w k (t)=0$ for $t<0$ and such that $\w K(\cdot)=\Z\w
k\in H^\infty$.
\end{definition}
\def\S{{\bf s}}
For $r\in[1,+\infty]$, let $\S:\ell_r\to\ell_r$ be the shift operator defined as $(\S x)(t)=x(t+1)$.
\begin{definition}\label{defP} Let  $\Y\subset \elr$ be a class of one-sided sequences, $r\in[1,+\infty]$.
\begin{itemize}
\item[(i)] We say that this class is unilaterally   $\elr$-predictable if
there exists a sequence $\{\w k_m(\cdot)\}_{m=1}^{+\infty}\subset
\w\K$ such that \baaa \|\S x-\w x_m\|_{\ell_r(-\infty,-1)}\to 0\quad
\hbox{as}\quad m\to+\infty\brea\quad\forall x\in\Y. \eaaa Here $  \w
x_m(t)\defi \sum^t_{s=-\infty}\w k_m(t-s)x(s).$
\item[(ii)]
 We say that the class $\Y$ is unilaterally  uniformly $\elr$-predictable  if, for any $\e>0$, there exists $\w k(\cdot)\in \w\K$ such that \baaa \|\S x- \w x\|_{\ell_r(-\infty,-1)}\le \e\quad
\forall x\in\Y. \label{predu}\eaaa Here $\w x(t)\defi
\sum^t_{s=-\infty}\w k(t-s)x(s).$
\end{itemize}
\end{definition}
\def\BL{{\scriptscriptstyle BL}}
\par
For $r\in\{1,2\}$, $\O\in[0,\O)$, let $\ell_{r,\BL}^-(\Omega)$ be the set of all one-sided sequences causally band-limited $x\in \elm$ such that, for each $x \in\ell_{r,\BL}^-(\Omega)$, the condition of  Definition \ref{defBL}
are satisfied. For $d>0$, let $\ell_{r,\BL}^-(\Omega,d)$ be the set of all $x\in \ell_{r,\BL}^-(\O)$ such that $\|x\|_{\elr}\le d$.
Let $\ell_{r,\BL}^-=\cup_{\O\in[0,\pi)}\ell_{r,\BL}^-(\Omega)$.
  \begin{theorem}\label{ThBL}
\begin{itemize}
\item[(i)]
The class $\ell_{2,\BL}^-$ is  unilaterally $\ell_2^-$-predictable.
\item[(ii)] For any  $\O\in[0,\pi)$, $d>0$, the class $\ell_{2,\BL}^-(\O,d)$  is unilaterally  uniformly
$\ell_2^-$-predictable.
\item[(iii)] Let  $\mu>1$  and $q>1$ be given.
A sequence of predicting kernels that ensures prediction required in
(i) and (ii) can be constructed as the following: $ \w k(\cdot)=\w
k(\cdot,\g)=\Z^{-1}\w K$, where \baa&&\w K(z)\breakk=z\left(
1-\exp\left[-\frac{\g}{z+ 1-\g^{2\mu/(1-q)} }\right]\right).\hphantom{xxx} \label{wK}\eaa
Here  $\g>0$ is a parameter; the
prediction error vanishes as $\g\to +\infty$.
\end{itemize}
\end{theorem}
The
predicting kernels (\ref{wK}) were suggested in \cite{D12a}, Theorem 1  for two-sided sequences.
They represent an extension
on the discrete time setting  of the construction introduced in \cite{D10} for continuous time processes.
\par
Further, let some $q>1$  be given. For $c>0$ and $\o\in[-\pi,\pi]$, set \baaa
h(\o,c)=\exp\frac{c}{[(\cos(\o)+1)^2+\sin^2(\o)]^{q/2}}.
\label{hdef}\eaaa
\par
Let $\W(c)$ be the class of all sequences $x\in\elm$ such that, for each $x\in\elm$,
at least one of the following conditions holds: either there exist $a\in\R$ such that
\baaa \esssup_{\o\in[0,\pi]} |\xi_1(\o)-a| h(\o,c)< +\infty,
\label{hfin}\eaaa
or
\baaa \esssup_{\o\in[0,\pi]} |\xi_2'(\o)| h(\o,c)< +\infty.
\label{hfin2}\eaaa
Here $\xi_k=\Xi_k x$, $\xi_2=(\xi_2'(\o),\xi_2'')$.
Let $\W=\cup_{c>0}\W(c)$. \par
 Note that $h(\o,c)\to +\infty$ as $\o\to\pm \pi$ and that, for $x\in \W(c)$, either $\xi_1(\o)-a$ or $\xi_2''(\o)$
vanishing with a sufficient rate of decay as $\o\to \pi$.
In particular, $\ell_{2,\BL}^-\subset\W$, i.e., the class $\W$ includes all causally band-limited sequences.
\par Further, for $c>0$  and $d>0$,  let $\V(c,d)$   be a class of processes $x\in \W(c)$ such that there exists $a\in\R$  such that
\baa  \min\Bigl(\esssup_{\o\in[0,\pi]} |(\xi_1(\o)-a)h(\o,c)|,\break\,\,\esssup_{\o\in[0,\pi]} |\xi''_2(\o)h(\o,c)|\Bigr)\le d, \label{V}\eaa where
$\xi_k=\Xi_k x$ and $\xi_2=(\xi_2'(\o),\xi_2'')$.
\par

%>> x=0.98; y= acos(-x); abs(exp(j*y)+1)^(-1)-(2-2*x)^(-1/2) ans = 1.7764e-015
  \begin{theorem}\label{ThM}
Let either $r=2$ or $r=+\infty$. \begin{itemize}
\item[(i)]
The class $\W$ is unilaterally  $\elr$-predictable.
\item[(ii)] For any given  $c>0$ be given and $d>0$, the class $\V(c,d)$ is uniformly unilaterally
$\ell_r$-predictable.
\item[(iii)]
A sequence of predicting kernels that ensures prediction required in
(i) and (ii) can be constructed  as defined in (\ref{wK}).
\end{itemize}
\end{theorem}
%\newpage
\section{Proofs}
{\em Proof of Lemma \ref{lemma1}}. For $x\in\elm$, we define $\ww x_k\in\ell_2$, $k=1,2$, such that $x_k(t)=x(t)$ for $t\le 0$ and
$\ww x_1(t)=x(-t)$, $\ww x_2(t)=-x(-t)$ for $t> 0$. For   $\xi_k=\Xi_k x$ and $X_k=\Z\ww x_k$, it can be verified that
\baa
X_1\ew=\xi_1(|\o|),\quad\o\in[-\pi,\pi],\label{XX1}\eaa
and, for $\xi_2=(\xi_2'(\o),\xi_2'')$,
 \baa
&&X_2\ew=i\xi_2'(\o)+\xi_2'',\quad  \o\in[0,\pi],\nonumber\\
&&X_2\ew=-i\xi_2'(-\o)+\xi_2'',\quad  \o\in[-\pi,0).\hphantom{xxx}\label{XX2}\eaa Since $\ww x_k\in\ell_2$, it follows that $|X_k\ew|\in L_2([-\pi,\pi],\R)$. Hence $\xi_k=\Xi_k x\in L_2([0,\pi],\R)$, and the mappings $\Xi_1:\elm\to L_2([0,\pi],\R)$ and $\Xi_2:\elm\to L_2([0,\pi],\R)\times\R$ are continuous.

Further, let $\xi_1\in L_2([0,\pi],\R)$ and $\xi_2=(\xi_2'(\o),\xi_2'')\in L_2([0,\pi],\R)\times\R$  be  arbitrarily selected.
Let us define mappings $X_k:\T\to \R$ defined by (\ref{XX1}) and (\ref{XX2}).
Let $\ww x_k=\Z^{-1}X_k$, i.e.,
\baaa\ww x_k(t)=\frac{1}{2\pi}\int_{-\pi}^\pi X_k\ew
e^{i\o t}d\o, \quad t=0,\pm 1,\pm 2,....\eaaa
It can be verified immediately that $\ww x_k(t)$ are real,
$\ww  x_1(t)=\ww x_1(-t)$,  $\ww x_2(t)=-\ww x_2(-t)$, $t>0$, and that (\ref{invinv})
holds for $x_k\defi \w x_k|_{\{t\le 0\}}$.
It follows that $\Xi_k x_k=\xi.$ Hence
the mappings $\Xi_1: \elm\to L_2([0,\pi],\R)$ and $\Xi_2: \elm\to L_2([0,\pi],\R)\times \R $ are bijections.
It is known that an inverse of a continuous bijection between Banach spaces is also continuous.
This completes the proof of  Lemma \ref{lemma1}. $\Box$
\par
{\em Proof of Theorem \ref{Th0}}.
It follows from the proof of Lemma \ref{lemma1} that if the condition of the theorem is satisfied that
one of the processes $\ww x_k$ introduced in this proof is band-limited.
This completes the proof of  Theorem  \ref{Th0}. $\Box$ \par
{\em Proof of Remark \ref{rem1}.}
It follows from the proof of Theorem \ref{Th0} that the  conditions on $\xi_1$ and $\xi_2'$ required in Theorem \ref{Th0} cannot be satisfied simultaneously; this would contradict the uniqueness of of the extrapolation of a band-limited processes established in \cite{D12a,D12b}.
\par
Theorem \ref{ThBL} follows immediately from Theorem \ref{ThM}. Therefore, it suffices to prove Theorem \ref{ThM}.
  \par
 {\em Proof of Theorem \ref{ThM}}. Let  $x\in\W$ and $\xi_k=\Xi_kx$, and let $\ww x_k$ and $ X_k$
 be such as defined in the proof of Lemma \ref{lemma1} above; in particular, we assume that (\ref{XX1})-(\ref{XX2}) hold.

 Let us prove statement (i). Let $x\in\W(c,d)$ be given for some $c>0$, $d>0$.
 If the definition of $\W(c,d)$ holds for this $x$
 with the condition for $\xi_1(\o)$, we set $k=1$ and select $a\in\R$ to be the corresponding $a$. Otherwise, we set $k=2$ and $a=x(0)$; this happens  if  the definition of $\W(c)$ holds for  $x$
 with the condition for $\xi''_2(\o)$ only.  Consider a process $x^a(t)$ such that $x^a(t)=x(t)$ for $t<0$ and $x^a(0)=x(0)-a$.
 Let $\xi^a_k=\Xi_k x^a$. By the definitions, it follows that
 $\esssup_{\o\in[0,\pi]} |\zeta^a_k(\o)| h(\o,c)< +\infty$, where $\zeta_1^a=\xi_1$ and
  $\zeta_2^a=\xi_2''$, for $\xi_2=(\xi_2'(\o),\xi_2'')$.
%\label{hfin2}\eaaa
\par
Let $\ww x_k^a(t)\in\ell_r$ be defined similarly to $\ww x_k$ in  in the proof of Lemma \ref{lemma1}, with $x$ replaced by $x^a$.
 By the definitions,  $\ww x^a_k$ belongs to the class $\X\subset\ell_r$  used in Theorem 1 (i)
 \cite{D12a}. By this theorem, the class $\X$ is  $\ell_r$-predicable
 in the sense of Definition 2(i) \cite{D12a}. Since $\ww x^a_k(t)=x(t)$ for $x\in\W$, $t< 0$, it follows that the predicability of $\X$ in the sense of Definition 2(i) \cite{D12a} implies predictability of $\W$ in the sense of Definition \ref{defP}(i).
\par
 Let us prove statement (ii) following the same steps.  Let $x\in\V(c,d)$. It suffices to consider $d=1$ only.
  If the definition of $\V(c,1)$ holds for this $x$
 with the minimum in (\ref{V}) achieved for $\xi_1(\o)$, we set $k=1$ and select $a\in\R$ to be the corresponding $a$. Otherwise, we set $k=2$ and $a=x(0)$; this happens  if  the  minimum in (\ref{V}) is achieved for $\xi''_2(\o)$ only.
  \par
  Let
 $x^a$, $\xi^a_k=\Xi_k x^a$, and $\ww x_k^a(t)\in\ell_r$, be defined as in the proof for the statement (i) above.

 By the definitions, $\ww x_k^a$ belong to the class $\U(c)$ defined in Theorem 1 (ii) \cite{D12a}. The class $\U(c)\subset\ell_r$ is uniformly $\ell_r$-predicable
 in the sense of Definition 2(ii) \cite{D12a}. Since $\ww x_k^a(t)=x(t)$ for $x\in\V(c,1)$, $t< 0$, it follows that the predicability for $\U(c_0)$
  in the sense of  Definition 2(ii) \cite{D12a}
implies predictability in the sense of Definition \ref{defP}(ii) for $\V(c,1)$. This completes the proof of  Theorem \ref{ThM}. $\Box$
\section{Discussion and future developments}
\begin{enumerate}
\item It is possible to consider complex valued sequences, with small modification of the definitions.
\item Similarly to the classical $Z$-transforms,  the transforms $\Xi_k$   allow to convert a linear difference equation in $\ell_r^-$
into an algebraic equation  and express the solution of this equations via the corresponding transfer function. It can be done via an appropriate extension of the equations on $t>0$ as backward difference equations.
\item A straightforward modification leads to analogs of transforms $\Xi_k$ for sequences from $\ell_r(0,+\infty)$.
\item
Theorems \ref{ThBL}-\ref{ThM} can be used for prediction of the sequences for times $t>0$ based on the observations for
$t\le 0$. These theorems allow to
 tell non-predicable  sequences $\{x(s)\}_{s=-\infty}^0$ from
predictable ones from $\ell_{2,\BL}^-$ or $\W$. For a given $\O\in[0,\pi)$, the set $\ell_{2,\BL}^-(\O)$
is closed in $\elm$.  The sequences from these set allows an extrapolation for $t>0$ such described in the proof of
Lemma \ref{lemma1} and such that
resulting two-sided sequence is band-limited. It could be natural
to consider the projection of an arbitrary   sequence on one of these subspaces, and accept the
extrapolation of this projection as the optimal forecast, following
the approach from \cite{D13,F,Mi,N,NT}.\index{in
Frank and Klotz (2002),  Miamee (1988), Nakazi (1984), Nakazi and
Takahashi (1980), and \cite{D13}.}
\item It is still an open question how far are away the sufficient conditions established
in Theorem \ref{Th0} from the necessary conditions of band-limitness.
\end{enumerate} % \end{document}
\index{\subsection*{Acknowledgment} This work  was
supported by ARC grant of Australia DP120100928 to the author.}

\end{document}